\documentclass{amsart}

\usepackage{amssymb} \usepackage{amsthm} \usepackage{amsmath}
\usepackage[all]{xy} \SelectTips{eu}{}
\usepackage[hidelinks]{hyperref}
\newcommand{\numberseries}{\bfseries} %Fontseries used for numbering

\newlength{\thmtopspace} %Space above theorem
\newlength{\thmbotspace} %Space below theorem
\newlength{\thmheadspace} %Space after theorem label
\newlength{\thmindent} %For indenting

\newtheoremstyle{bfupright head,upright body}
{\thmtopspace}{\thmbotspace}
{\upshape}{\thmindent}{\bfseries}{.}{\thmheadspace} {{\numberseries
    \thmnumber{#2\;}}\thmnote{#3}}

\newtheoremstyle{fixed bf head,slanted body}
{\thmtopspace}{\thmbotspace}{\slshape}
{\thmindent}{\bfseries}{.}{\thmheadspace} {{\numberseries
    \thmnumber{#2\;}}\thmname{#1}\thmnote{ (#3)}}

\newtheoremstyle{fixed bf head,upright body}
{\thmtopspace}{\thmbotspace}{\upshape}
{\thmindent}{\bfseries}{.}{\thmheadspace} {{\numberseries
    \thmnumber{#2\;}}\thmname{#1}\thmnote{ (#3)}}

\theoremstyle{bfupright head,upright body} \newtheorem{res}{}[section]

\theoremstyle{fixed bf head,slanted body}
\newtheorem{thm}[res]{Theorem} \newtheorem*{thm*}{Theorem}
\newtheorem{prp}[res]{Proposition} \newtheorem*{prp*}{Proposition}
\newtheorem{cor}[res]{Corollary} \newtheorem*{cor*}{Corollary}
 \newtheorem*{lem*}{Lemma}

\theoremstyle{fixed bf head,upright body}
 \newtheorem*{dfn*}{Definition}
\newtheorem{rmk}[res]{Remark} \newtheorem*{rmk*}{Remark}
\newtheorem{exa}[res]{Example}
\newtheorem{qst}[res]{Question} 

\newlength{\thmlistleft} %leftmargin
\newlength{\thmlistright} %rightmargin
\newlength{\thmlistpartopsep} %partopsep
\newlength{\thmlisttopsep} %topsep
\newlength{\thmlistparsep} %parsep
\newlength{\thmlistitemsep} %itemsep

\newcounter{eqc} \newenvironment{eqc}{\begin{list}{\upshape
      (\textit{\roman{eqc}})}%
    {\usecounter{eqc}%
      \setlength{\leftmargin}{\thmlistleft}%
      \setlength{\labelwidth}{\thmlistleft}%
      \setlength{\rightmargin}{\thmlistright}%
      \setlength{\partopsep}{\thmlistpartopsep}%
      \setlength{\topsep}{\thmlisttopsep}%
      \setlength{\parsep}{\thmlistparsep}%
      \setlength{\itemsep}{\thmlistitemsep}}}%
  {\end{list}}%

\newcommand{\eqclbl}[1]{{\upshape(\textit{#1})}}

\newenvironment{prf*}[1][Proof]{%
  \begin{proof}[\bf #1]
    \setcounter{equation}{0}
    } {\end{proof} }

\newcommand{\proofofimp}[3][:]{\mbox{\eqclbl{#2}$\!\implies\!$\eqclbl{#3}#1}}

\newcommand{\thmref}[2][Theorem~]{#1\ref{thm:#2}}
\newcommand{\corref}[2][Corollary~]{#1\ref{cor:#2}}
\newcommand{\prpref}[2][Proposition~]{#1\ref{prp:#2}}

\newcommand{\rmkref}[2][Remark~]{#1\ref{rmk:#2}}

\newcommand{\thmcite}[2][?]{\cite[Theorem~#1]{#2}}
\newcommand{\corcite}[2][?]{\cite[Corollary~#1]{#2}}
\newcommand{\prpcite}[2][?]{\cite[Proposition~#1]{#2}}
\newcommand{\lemcite}[2][?]{\cite[Lemma~#1]{#2}}
\newcommand{\rmkcite}[2][?]{\cite[Remark~#1]{#2}}
\newcommand{\dfncite}[2][?]{\cite[Definition~#1]{#2}}
\newcommand{\setof}[3][\mspace{1mu}]{\{#1#2 \mid #3#1\}}
\newcommand{\ZZ}{\mathbb{Z}} \newcommand{\NN}{\mathbb{N}}
\newcommand{\deq}{\:=\:} 
\newcommand{\dis}{\:\is\:}
\renewcommand{\d}{v} % generic homological degree
\newcommand{\is}{\cong} 
\renewcommand{\le}{\leqslant} \renewcommand{\ge}{\geqslant}
\newcommand{\lra}{\longrightarrow}

 \newcommand{\Aop}{A^\circ}

\newcommand{\Cy}[2][]{\operatorname{Z}_{#1}(#2)}
\newcommand{\Co}[2][]{\operatorname{C}_{#1}(#2)}
\renewcommand{\H}[2][]{\operatorname{H}_{#1}(#2)}

\newcommand{\fd}[2][A]{\operatorname{fd}_{#1}#2}
\newcommand{\id}[2][A]{\operatorname{id}_{#1}#2}

\newcommand{\Gfd}[2][A]{\operatorname{Gfd}_{#1}#2}
\newcommand{\Gid}[2][A]{\operatorname{Gid}_{#1}#2}
\newcommand{\Gpd}[2][A]{\operatorname{Gpd}_{#1}#2}
\newcommand{\Gfcd}[2][A]{\operatorname{Gfcd}_{#1}#2}
\newcommand{\Ggldim}[1][A]{\operatorname{Ggldim}(#1)}
\newcommand{\Gwgldim}[1][A]{\operatorname{Gwgldim}(#1)}
\newcommand{\FPD}[1][A]{\operatorname{FPD}(#1)}
\newcommand{\FFD}[1][A]{\operatorname{FFD}(#1)}

\newcommand{\sfli}[1][A]{\operatorname{sfli}(#1)}
\newcommand{\Hom}[3][A]{\operatorname{Hom}_{#1}(#2,#3)}
\newcommand{\Ext}[4][A]{\operatorname{Ext}_{#1}^{#2}(#3,#4)}
\newcommand{\Ker}[1]{\operatorname{Ker}#1}
\newcommand{\Coker}[1]{\operatorname{Coker}#1}
\def\urltilda{\kern -.15em\lower .7ex\hbox{\~{}}\kern .04em}

\title{Five theorems on Gorenstein global
  dimensions} %Rings of Gorenstein (weak) global dimension zero

\author[L.W.\ Christensen]{Lars Winther Christensen} %
\address{L.W.C. \ Texas Tech University, Lubbock, TX 79409, U.S.A.}
\email{lars.w.christensen@ttu.edu}
\urladdr{http://www.math.ttu.edu/\urltilda lchriste}

\author[S.\ Estrada]{Sergio Estrada} %
\address{S.E. \ Universidad de Murcia, Murcia 30100, Spain}
\email{sestrada@um.es} %
\urladdr{https://webs.um.es/sestrada/}

\author[P.\ Thompson]{Peder Thompson} %
\address{P.T. \ Niagara University, Niagara University, NY 14109,
  U.S.A.}  \email{thompson@niagara.edu} %
\urladdr{http://pthompson.nupurple.net}

\thanks{L.W.C.\ was partly supported by Simons Foundation
  collaboration grant 428308. S.E. was partly supported by grant
  PID2020-113206GB-I00 funded by MCIN/AEI/10.13039/ 501100011033.}

\date{10 August 2022}

\keywords{Ding-Chen ring; Gorenstein global dimension; Gorenstein weak global
  dimension; Gorenstein flat-cotorsion dimension; IF ring;
  Quasi-Frobenius ring; strongly cotorsion module}

\subjclass[2020]{Primary 16E10. Secondary 16E65; 16S34.}

\begin{document}

\begin{abstract}
  We expand on two existing characterizations of rings of Gorenstein
  (weak) global dimension zero and give two new characterizations of
  rings of finite Gorenstein (weak) global dimension.  We also include
  the answer to a question of Y.~Xiang on Gorenstein weak global
  dimension of group rings.
\end{abstract}

\maketitle

%%%%%%%%%%%%%%%%%%%%%%%%%%%%%%%%%%%%%%%%%%%%%%%%%%%%%%%%%%%%%%%%%%%%%%

\section*{Introduction}

\noindent
In this paper, $A$ denotes a unital associative ring. By an $A$-module
we mean a left $A$-module; right $A$-modules are considered as modules
over the opposite ring $\Aop$.

Two classic facts in algebra are that $A$ has finite weak global
dimension if and only if so does $\Aop$, and that $A$ has global
dimension zero if and only if so does $\Aop$. Corresponding results
are known to hold in the Gorenstein setting: In \cite{CET-21a} we
observe that the Gorenstein weak global dimension is symmetric, and in
\cite{DBnNMh10} Bennis and Mahdou show a ring has Gorenstein global
dimension zero if and only if it is Quasi-Frobenius, which is a
left--right symmetric property.

Gorenstein flat-cotorsion modules, defined in \cite{CET-20}, and the
associated dimension developed in \cite{CELTWY-21} turn out to be
useful for refining our understanding of Gorenstein global
dimensions. Two of our main results, \thmref[Theorems~]{IF} and
\thmref[]{QF}, characterize rings of Gorenstein (weak) global
dimension zero. In the first, we show that IF rings are precisely
those rings where all modules have Gorenstein flat-cotorsion dimension
at most zero. In the second, we show that Quasi-Frobenius rings are
those rings where all modules are Gorenstein flat-cotorsion. Further,
we show in \thmref[Theorems~]{Ggl_SC_GI} and \thmref[]{Gfcd_SC_GI}
that under assumptions of finite finitistic dimensions, rings of
finite Gorenstein (weak) global dimension are precisely the rings over
which strongly cotorsion modules and Gorenstein injective modules
coincide.

An important ingredient in the proofs of these results are
generalizations of two results due to Holm \cite{HHl04c}: We prove in
\thmcite[2.1]{CET-21a} that the flat, Gorenstein flat, and Gorenstein
flat-cotorsion dimensions agree for modules of finite injective
dimension, and in \thmref{idGid} we show that the injective and
Gorenstein injective dimensions agree for modules of finite flat
dimension.
%\enlargethispage*{\baselineskip}
In the final section, \thmref{grouprings} provides an answer to a
question of Y.~Xiang \cite{YXn21} on the Gorenstein global dimension
of group rings.

% \begin{equation*}
%   \ast \ \ast \ \ast
% \end{equation*}
% \noindent
The notation employed in the paper is standard: we write pd, id, and
fd for the projective, injective, and flat dimensions of modules and
complexes. The corresponding Gorenstein dimensions are denoted Gpd,
Gid, and Gfd. The Gorenstein flat-cotorsion dimension from
\cite{CELTWY-21} gets abbreviated Gfcd. We write $\FPD$ and $\FFD$ for
the finitistic projective and finitistic flat dimension of $A$. The
invariant $\sfli$ is defined as 
$\sup\{\fd{M}\mid M \text{ is an injective $A$-module}\}$. The
Gorenstein global dimension is defined as 
\begin{equation*}
  \Ggldim = \sup\{\Gpd{M} \mid M \text{ is an $A$-module}\} \:,
\end{equation*}
and it equals $\sup\{\Gid{M} \mid M \text{ is an $A$-module}\}$. The
Gorenstein weak global dimension is defined similarly,
\begin{equation*}
  \Gwgldim = \sup\{\Gfd{M} \mid M \text{ is an $A$-module}\} \:.
\end{equation*}

In dealing with $A$-complexes we use homological notation. For an
$A$-complex $M$ and $n \in \ZZ$ set $\Cy[n]{M} = \Ker{\partial_n^M}$
and $\Co[n]{M} = \Coker{\partial^M_{n+1}}$ so that they are,
respectively, submodules and quotient modules of the module in degree
$n$.

%%%%%%%%%%%%%%%%%%%%%%%%%%%%%%%%%%%%%%%%%%%%%%%%%%%%%%%%%%%%%%%%%%%%%%
\section{Gorenstein injective dimension of flat modules}
\label{sec_Gid}

\noindent
Recall that an $A$-module $M$ is \emph{cotorsion} if $\Ext{1}{F}{M}=0$
holds for every flat $A$-module $F$, and $M$ is \emph{strongly
  cotorsion} if $ \Ext{1}{L}{M}=0$ holds for every $A$-module $L$ of
finite flat dimension; see Xu \dfncite[5.4.1]{xu}.

By convention the homological supremum of an acyclic complex is
$-\infty$, so the next result says, in particular, that every cycle
module in an acyclic complex of strongly cotorsion modules is strongly
cotorsion.

\begin{prp}
  \label{prp:GI_is_SC}
  Let $M$ be a complex of strongly cotorsion $A$-modules. For every
  $n \ge \sup\{v \in \ZZ \mid \H[\d]{M}\ne 0 \}$ the module
  $\Co[n]{M}$ is strongly cotorsion.
\end{prp}

\begin{prf*}
  If $M$ is acyclic, then $\Co[n]{M}$ for every $n\in\ZZ$ is a
  cokernel in an acyclic complex of strongly cotorsion modules; we
  first reduce the general case to this special case. Set
  $s = \sup\{v \in \ZZ \mid \H[\d]{M}\ne 0 \}$. Splicing a shifted
  injective resolution of the module $\Co[s]{M}$ with the acyclic
  complex $\cdots \to M_{s+1} \to M_s \to \Co[s]{M} \to 0$ one gets an
  acyclic complex $I$ with $\Co[n]{I} = \Co[n]{M}$ for $n \ge s$. It
  now suffices to show that the modules $\Co[i]{I}$ are strongly
  cotorsion.

  Let $L$ be a module with $\fd{L}\le f$. There is an exact sequence,
  \begin{equation*} 0 \lra F_f \lra \cdots \lra F_0 \lra L \lra 0 \:,
  \end{equation*}
  with each $F_i$ a flat $A$-module. Given a cotorsion module $C$,
  dimension shifting along this sequence yields
  $\Ext{f+1}{L}{C} \cong \Ext{1}{F_f}{C}=0$. By a result of Bazzoni,
  Cort\'es-Izurdiaga, and Estrada \thmcite[1.3]{BCE-20} each module
  $\Co[i]{I}$ is cotorsion, so for every $i \in \ZZ$ one has
  $\Ext{f+1}{L}{\Co[i]{I}}=0$. Now, for every $i\in\ZZ$ there is an
  exact sequence,
  \begin{equation*}
    0 \lra \Co[i+f]{I} \lra I_{i+f-1} \lra \cdots \lra I_{i} \lra \Co[i]{I} \lra 0 \:.
  \end{equation*}
  Dimension shifting yields
  $\Ext{1}{L}{\Co[i]{I}}\cong \Ext{f+1}{L}{\Co[i+f]{I}}=0$, whence
  $\Co[i]{I}$ is strongly cotorsion.
\end{prf*}

Recall that an acyclic complex $I$ of injective $A$-modules is called
\emph{totally acyclic} if $\Hom{E}{I}$ is acyclic for every injective
$A$-module $E$. Recall further that an $A$-module $G$ is called
\emph{Gorenstein injective} if there exists a totally acyclic complex
$I$ of injective $A$-modules with $\Cy[0]{I} \is G$. It follows from
\prpref{GI_is_SC} that a Gorenstein injective $A$-module is strongly
cotorsion.

Now, if a flat $A$-module $F$ is Gorenstein injective, then one gets
from the defining totally acyclic complex an exact sequence
$0 \to G \to I \to F \to 0,$ with $I$ injective and $G$ Gorenstein
injective. Since $\Ext{1}{F}{G}=0$ holds this sequence splits, which
means that $F$ is a summand of $I$ and hence injective. This argument
can be souped up to yield the next theorem, which is dual to
\thmcite[1.1]{CET-21a} and subsumes \thmcite[2.1]{HHl04c}

Recall, for example from \prpcite[3.6]{CKL-17}, that the Gorenstein
injective dimension of an $A$-complex $M$ can be defined as the least
integer $n$ such that $(1)$ $\H[\d]{M} = 0$ holds for all $\d < -n$
and $(2)$ there exists a semi-injective $A$-complex $I$, isomorphic to
$M$ in the derived category, such that the cycle submodule
$\Cy[-n]{I}$ is Gorenstein injective.

\begin{thm}
  \label{thm:idGid}
  Let $M$ be an $A$-complex with bounded homology. If\, $\fd{M}$ is
  finite, then $\Gid{M} = \id{M}$ holds.
\end{thm}

\begin{prf*}
  The equality $\Gid{M} = \id{M}$ holds trivially if $M$ is acyclic,
  so assume that $M$ is not acyclic and, without loss of generality,
  assume further that $\fd{M}=0$ holds.  Set
  $u = \min\setof{v\in\ZZ}{\H[\d]{M} \ne 0}$ and let $F$ be a bounded
  complex of flat $A$-modules with $F_v=0$ for $v > 0$ and $v < u$,
  such that $F$ and $M$ are isomorphic in the derived category of $A$.
  As one has $\Gid{F} \le \id{F}$ it suffices to show that
  $\Gid{F} =n$ implies $\id{F} \le n$. Let $F \to I$ be a
  semi-injective resolution with $I_v=0$ for $v>0$ and $C$ be its
  mapping cone. As $I_v=0$ holds for $v>0$ one has $C_1 = F_0$ and
  $C_v=0$ for $v>1$. As $F_v=0$ holds for $v<u$ one has
  \begin{equation}
    \tag{1}
    \Cy[v]{C} = \Cy[v]{I}  \quad\text{for}\quad v \le u\:.
  \end{equation}
  Assume that $\Gid{F} = n$ holds; one then has $-n \le u$. Let $G$ be
  a Gorenstein injective $A$-module. From the defining totally acyclic
  complex of injective $A$-modules one gets an exact sequence,
  \begin{equation*}
    0 \lra H \lra I_0 \lra \cdots \lra I_{-n} \lra G \lra 0 \:,
  \end{equation*}
  with $H$ a Gorenstein injective $A$-module. Dimension shifting along
  this exact sequence yields
  \begin{equation}
    \tag{2}
    \Ext{1}{\Cy[-(n+1)]{C}}{G} \dis \Ext{n+2}{\Cy[-(n+1)]{C}}{H} \:.
  \end{equation}
  Since $H$ is Gorenstein injective, and hence cotorsion by
  \prpref{GI_is_SC}, and the modules $C_i$ are direct sums of flat
  modules and injective modules, dimension shifting along the exact
  sequence
  \begin{equation*}
    0 \lra F_0 \lra C_0 \lra \cdots \lra C_{-n} \lra \Cy[-(n+1)]{C} \lra 0
  \end{equation*}
  yields
  \begin{equation}
    \tag{3}
    \Ext{n+2}{\Cy[-(n+1)]{C}}{H} \dis \Ext{1}{F_0}{H} \deq 0 \:.
  \end{equation}
  Combining $(1)$--$(3)$ one gets $\Ext{1}{\Cy[-(n+1)]{I}}{G} =0$ for
  every Gorenstein injective $A$-module $G$. In particular
  $\Ext{1}{\Cy[-(n+1)]{I}}{\Cy[-n]{I}} =0$ holds, so the exact
  sequence $0 \to \Cy[-n]{I} \to I_{-n} \to \Cy[-(n+1)]{I} \to 0$
  splits, which means that $\Cy[-n]{I}$ is injective and
  $\id{F} \le n$ holds as desired.
\end{prf*}

%%%%%%%%%%%%%%%%%%%%%%%%%%%%%%%%%%%%%%%%%%%%%%%%%%%%%%%%%%%%%%%%%%%%%%
\section{IF rings}
\label{sec_IFrings}

\noindent Following Colby~\cite{RRC75} a ring over which every
injective module is flat is called a left IF ring. If a ring $A$ and
its opposite ring $\Aop$ are both left IF, then $A$ is called
IF. Bennis \cite{DBn10} characterized IF rings in terms of the
Gorenstein weak global dimension, which at that point was not known to
be symmetric. In the commutative case, where this distinction is
irrelevant, the characterization was also obtained by Mahdou,
Tamekkante, and Yassemi~\cite{MTY-11}. Here we use the Gorenstein
flat-cotorsion dimension to characterize left-IF rings; Bennis'
characterization \prpcite[2.14]{DBn10} is recovered as the equivalence of \eqclbl{i} and \eqclbl{iii} in \thmref{IF}.

Recall that an acyclic complex $F$ of flat-cotorsion $A$-modules is
called \emph{totally acyclic} if $\Hom{F}{C}$ is acyclic for every
flat-cotorsion $A$-module $C$. Recall further that an $A$-module $G$
is called \emph{Gorenstein flat-cotorsion} if there exists a totally
acyclic complex $F$ of flat-cotorsion $A$-modules with
$\Co[0]{F} \is G$. It follows from \thmcite[1.3]{BCE-20} that a
Gorenstein flat-cotorsion $A$-module is cotorsion.

Recall from \dfncite[4.1]{CELTWY-21}, that the Gorenstein
flat-cotorsion dimension of an $A$-complex $M$ can be defined as the
least integer $n$ such that $(1)$ $\H[\d]{M} = 0$ holds for all
$\d > n$ and $(2)$ there exists a semi-flat-cotorsion $A$-complex
$F$---i.e.\ a semi-flat complex consisting of flat--cotorsion
modules---isomorphic to $M$ in the derived category, such that the
cokernel $\Co[n]{F}$ is Gorenstein flat-cotorsion.

\begin{prp}
  \label{prp:leftIF}
  The following conditions are equivalent.
  \begin{eqc}
  \item $\Gfcd{M} \le 0$ holds for every $A$-module $M$.
  \item An $A$-module is flat-cotorsion if and only if it is
    injective.
  \item An $A$-module is cotorsion if and only if it is Gorenstein
    flat-cotorsion.
  \item An $A$-module is cotorsion if and only if it is Gorenstein
    injective.
  \item $A$ is a left IF ring and an $A$-module is Gorenstein
    flat-cotorsion if and only if it is Gorenstein injective.
  \end{eqc}
\end{prp}

\begin{prf*}
  \proofofimp{i}{ii} If $F$ is flat-cotorsion, then $\Ext{1}{M}{F}=0$
  holds for every $A$-module $M$ by \thmcite[4.5]{CELTWY-21}, which
  implies that $F$ is injective. Let $I$ be injective;
  \thmcite[1.1]{CET-21a} yields $\fd{I} = \Gfcd{I} \le 0$, so $I$ is
  flat-cotorsion.

  \proofofimp{ii}{iii} Every Gorenstein flat-cotorsion module is
  cotorsion. To prove the converse, let $C$ be a cotorsion
  $A$-module. It is a cycle submodule in an acyclic complex of
  flat-cotorsion $A$-modules: Indeed, take an injective resolution to
  the right and a flat-cotorsion resolution to the left.  By assumption,
  every injective module is flat-cotorsion, so this is an acyclic
  complex of flat-cotorsion modules. Moreover, every flat-cotorsion
  module is injective, so it is in fact a totally acyclic complex of
  flat-cotorsion modules; in particular, $C$ is Gorenstein
  flat-cotorsion.

  \proofofimp{iii}{i} Let $M$ be an $A$-module. There exists a
  semi-flat-cotorsion complex $F$ isomorphic to $M$ in the derived
  category; see Nakamura and Thompson \thmcite[A.6]{TNkPTh20}. Since
  the cokernel $\Co[0]{F}$ has a left resolution by cotorsion modules,
  it is itself cotorsion by \lemcite[5.6]{CELTWY-21}. Thus $\Co[0]{F}$
  is Gorenstein flat-cotorsion by assumption. By the definition of
  Gorenstein flat-cotorsion dimension it follows that $\Gfcd{M}\le 0$
  holds.

  \proofofimp{ii}{iv} By \prpref{GI_is_SC} every Gorenstein injective
  $A$-module is cotorsion. Now, let $C$ be cotorsion. A left
  resolution of $C$ constructed by taking flat covers consists of
  flat-cotorsion modules. Splice this resolution together with a right
  injective resolution of $C$. As flat-cotorsion modules are
  injective, this produces an acyclic complex of injective modules. In
  particular, the complex has cotorsion cycles, see
  \thmcite[1.3]{BCE-20}. As injective modules are flat, this complex
  is $\Hom{I}{-}$-exact for every injective $A$-module $I$. Thus $C$
  is Gorenstein injective.

  \proofofimp{iv}{ii} Let $I$ be an injective $A$-module. For every
  cotorsion module $C$, one has $\Ext{1}{I}{C}=0$ because $C$ is also
  Gorenstein injective. Thus $I$ is flat-cotorsion. Conversely, let
  $F$ be flat-cotorsion. By assumption, $F$ must then be Gorenstein
  injective, so $F$ is injective by \thmref{idGid}.

  \proofofimp{ii}{v} As injective $A$-modules are flat, $A$ is a left
  IF ring. Further, as flat-cotorsion modules are precisely the
  injective modules, Gorenstein flat-cotorsion modules and Gorenstein
  injective modules coincide.

  \proofofimp{v}{ii} Every injective module is flat and hence
  flat-cotorsion. A flat-cotorsion module $F$ is Gorenstein
  flat-cotorsion and hence Gorenstein injective, so \thmref{idGid}
  yields $\id{F} = \Gid{F} = 0$.
\end{prf*}

\begin{thm}
  \label{thm:IF}
  The following conditions are equivalent.
  \begin{eqc}
  \item $\Gwgldim = 0$.
  \item $\Gfcd{M} \le 0$ holds for every $A$-module $M$ and
    $\Gfcd[\Aop]{M} \le 0$ holds for every $\Aop$-module $M$.
  \item $A$ is an IF ring.
  \end{eqc}
\end{thm}

\begin{prf*}
  The Gorenstein weak global dimension is symmetric, see
  \corcite[2.5]{CET-21a}, so $(i)$ implies $(ii)$ in view of
  \thmcite[5.7]{CELTWY-21}. Further $(ii)$ implies $(iii)$ by
  \prpref{leftIF}. Finally, to see that $(iii)$ implies $(i)$ note
  that since every injective $\Aop$-module is flat, every acyclic
  complex of flat $A$-modules is \textbf{F}-totally acyclic, i.e.\ it
  remains exact when tensored by an injective $\Aop$-module. Let $M$ be
  an $A$-module. It suffices to build an acyclic complex $F$ of flat
  $A$-modules such that $\Cy[0]{F}=M$. The left half is built by taking
  successive flat covers, whereas the right half is obtained by taking
  an injective resolution of $M$ and observing that it is a complex of
  flat modules by assumption.
\end{prf*}

\begin{rmk}
  \label{rmk:gfcd}
  While condition $(i)$ in \thmref{IF} can be interpreted as saying
  that all $A$- and $\Aop$-modules are Gorenstein flat, $(ii)$ does
  not say that all $A$- and $\Aop$-modules are Gorenstein
  flat-cotorsion. In fact, an $A$-module $M$ with $\Gfcd{M}=0$ is
  Gorenstein flat-cotorsion if and only if it is cotorsion; see
  \rmkcite[4.6]{CELTWY-21}.
\end{rmk}

Colby \cite[Proposition 5]{RRC75} shows that a ring $A$ is von Neumann
regular if and only if it is left IF with
$\operatorname{wgldim}(A) <\infty$. Here is the Gorenstein analog of
this result:
  
\begin{cor}
  \label{cor:IF ring}
  The following conditions are equivalent.
  \begin{eqc}
  \item $A$ is an IF ring.
  \item $A$ is a left or right IF ring with $\Gwgldim < \infty$.
  \end{eqc}
\end{cor}

\begin{prf*}
  By \thmref{IF} and \corcite[2.5]{CET-21a} it suffices to show that a
  left IF ring $A$ of finite Gorenstein weak global dimension has
  $\Gwgldim = 0$. Let $A$ be such a ring and $M$ an $A$-module. An
  injective resolution of $M$ is a right resolution by flat
  $A$-modules, so for every $n$ the module $M$ is a flat syzygy of its
  $n^{\rm th}$ injective cosyzygy. Thus $M$ is Gorenstein flat.
\end{prf*}

\begin{rmk}
  There exist right IF rings which are not IF, see \cite[Example
  2]{RRC75}. Per \corref{IF ring} such rings must be of infinite
  Gorenstein weak global dimension.
\end{rmk}

%%%%%%%%%%%%%%%%%%%%%%%%%%%%%%%%%%%%%%%%%%%%%%%%%%%%%%%%%%%%%%%%%%%%%%
\section{Quasi-Frobenius rings}
\label{sec_QFrings}

\noindent
The difference between a module being Gorenstein flat-cotorsion and
having Gorenstein flat-cotorsion dimension zero was already commented
on in \rmkref{gfcd}. The former quality is the stronger one, and the
main result of this section is that Gorenstein flat-cotorsionness of
all modules characterizes rings of Gorenstein global dimension zero.
First we characterize left perfect rings in terms of Gorenstein
flat-cotorsion modules.

\begin{prp}
  \label{prp:perfect}
  The following conditions are equivalent.
  \begin{eqc}
  \item Every Gorenstein flat $A$-module is Gorenstein projective.
  \item Every flat $A$-module is Gorenstein projective.
  \item Every Gorenstein projective $A$-module is Gorenstein
    flat-cotorsion.
  \item Every projective $A$-module is cotorsion.
  \item $A$ is left perfect.
  \end{eqc}
  Moreover, when these conditions are satisfied, an $A$-module is
  Gorenstein projective if and only if it is Gorenstein
  flat-cotorsion.
\end{prp}

\begin{prf*}
  Notice first that if $A$ is left perfect, then every $A$-module is
  cotorsion and the definitions of Gorenstein projective $A$-modules
  and Gorenstein flat-cotorsion $A$-modules coincide. This justifies
  the last assertion as well as the implication $(v) \implies
  (iii)$. The implications $(i) \implies (ii)$ and
  $(iii) \implies (iv)$ are trivial.

  \proofofimp{ii}{iv} As projective $A$-modules are right
  Ext-orthogonal to Gorenstein projective modules, it follows that
  each projective $A$-module is right Ext-orthogonal to flat
  $A$-modules and hence cotorsion.

  \proofofimp{iv}{v} The free $A$-module $A^{(\NN)}$ is in particular
  cotorsion, whence $A$ is left perfect by Guil~Asensio and Herzog
  \corcite[20]{PGAIHz05}.

  \proofofimp{v}{i} Every $A$-module is cotorsion, so per
  \thmcite[5.2]{CELTWY-21} every Gorenstein flat $A$-module is
  Gorenstein flat-cotorsion, which over a left perfect ring means
  Gorenstein projective.
\end{prf*}

The equivalence of $(i)$, $(i')$, and $(iii)$ in the next result was
proved by Bennis and Mahdou in \cite[Proposition
2.6]{DBnNMh10}. Recall that a Quasi-Frobenius ring is one where the
projective modules are injective. Over such a ring, the projective and
injective modules coincide, and the same is true for the opposite
ring.

\begin{thm}
  \label{thm:QF}
  The following conditions are equivalent.
  \begin{eqc}
  \item $\Ggldim = 0$.
  \item[$(i')$] $\Ggldim[\Aop] = 0$.
  \item Every $A$-module is Gorenstein flat-cotorsion.
  \item[$(ii')$] Every $\Aop$-module is Gorenstein flat-cotorsion.
  \item $A$ is Quasi-Frobenius.
  \end{eqc}
\end{thm}

\begin{prf*}
  As the Quasi-Frobenius property is left--right symmetric, it
  suffices to show the equivalence of the unprimed conditions.

  \proofofimp{i}{ii} By assumption every $A$-module is Gorenstein
  projective, so it follows from \prpref{perfect} that every
  $A$-module is Gorenstein flat-cotorsion.

  \proofofimp{ii}{iii} Let $P$ be a projective $A$-module. As every
  module is Gorenstein flat-cotorsion, $P$ is flat-cotorsion and
  $\Ext{1}{M}{P} = 0$ holds for every $A$-module $M$, so $P$ is
  injective.

  \proofofimp{iii}{i} A Quasi-Frobenius ring is noetherian and
  self-injective on both sides, so it is in particular an
  Iwanaga-Gorenstein ring. Thus for every $A$-module $M$ one has
  $\Gpd{M} \le \id{A} =0$, see for example Holm \thmcite[2.28]{HHl04a}
  and Bass \prpcite[4.3]{HBS62}.
\end{prf*}

\section{Strongly cotorsion modules} %%%%%%%%%%%%%%%%%%%%%%%%%%%%%%%%%
\label{sec:SC}

\noindent
As noticed it follows from \prpref{GI_is_SC} that every Gorenstein
injective module is strongly cotorsion. A question considered in the
literature is: For which rings do Gorenstein injective modules and
strongly cotorsion modules coincide?

Yoshizawa \thmcite[2.7]{TYs12} shows that over (commutative)
Gorenstein complete local rings, Gorenstein injective modules and
strongly cotorsion modules are the same. Huang \thmcite[3.10]{ZHn13}
proves that the same statement holds for every Iwanaga-Gorenstein ring
$A$ such that the injective envelope of $A$ is flat, and Iacob
\thmcite[10]{AIc17} proves it for all Iwanaga-Gorenstein rings;
cf.~\prpcite[9.1.2]{rha}. The best result to date can be pieced
together from work of Gillespie:

\begin{rmk}
  \label{rmk:Gil}
  A ring is \emph{Ding-Chen} if it is coherent and has finite self
  FP-injective dimension on both sides.  Over a Ding-Chen ring, the
  classes of Gorenstein injective and strongly cotorsion modules
  coincide: Indeed, let $A$ be a Ding-Chen ring. Gillespie shows in
  \cite[Theorem 4.2 and Corollary 4.5]{JGl10} that the modules that
  are right $\operatorname{Ext}_A^1$-orthogonal to the class of
  $A$-modules with finite flat dimension---that is, the strongly
  cotorsion $A$-modules---coincide with the so-called Ding injective
  $A$-modules. Further, Gillespie shows \thmcite[1.1]{JGl17b} that
  Ding injective $A$-modules are the same as Gorenstein injective
  $A$-modules.
\end{rmk}

\begin{prp}
  \label{prp:supGfcd}
  If the quantity
  \begin{equation*}
    \sup\{\Gfcd M \mid M\text{ is an $A$-module}\}
  \end{equation*}
  is finite, then $\FFD$ is finite and an $A$-module is strongly
  cotorsion if and only if it is Gorenstein injective.
\end{prp}

\begin{prf*}
  Set $n= \sup\{\Gfcd M \mid M\text{ is an $A$-module}\}$ and assume
  that it is finite. This assumption implies that both $\FFD$ and
  $\sfli$ are finite: Indeed \thmcite[5.12]{CELTWY-21} yields
  $\FFD\le n$, and \thmcite[2.1]{CET-21a} yields $\sfli \le n$.

  As already noticed, every Gorenstein injective $A$-module is
  strongly cotorsion. Now let $M$ be a strongly cotorsion
  $A$-module. Since $\sfli$ is finite, it follows from
  \prpref{GI_is_SC} that every acyclic complex of injectives is
  totally acyclic. Therefore, it suffices to show that $M$ is the
  homomorphic image of an injective $A$-module with strongly cotorsion
  kernel. Since $\FFD$ is finite, results of Trlifaj \cite[Lemma
  1.5(3) and Theorem 1.14]{JTr07} apply to yield a short exact
  sequence,
  \begin{equation*}
    0\lra K\lra E\lra M\lra 0\;,
  \end{equation*} 
  with $K$ a strongly cotorsion $A$-module and $E$ an $A$-module of
  finite flat dimension. By the assumption on $M$ it follows that also
  $E$ is strongly cotorsion. To complete the proof we show that $E$ is
  an injective $A$-module. Consider an exact sequence,
  \begin{equation*}
    0 \lra E \lra I \lra C \lra 0 \;,
  \end{equation*}
  with $I$ an injective $A$-module. Since $E$ and $I$ have finite flat
  dimension, the $A$-module $C$ has finite flat dimension. But then,
  since $E$ is a strongly cotorsion $A$-module, the sequence splits
  and thus $E$ is injective.
\end{prf*}

A noetherian ring is Ding-Chen if and only if it is
Iwanaga-Gorenstein. Thus for noetherian rings, the next result is
equivalent to the statement in \rmkref{Gil}.

\begin{cor}
  \label{cor:Ggldim_fin_SCGI}
  If $\Ggldim$ is finite, then an $A$-module is strongly cotorsion if
  and only if it is Gorenstein injective.
\end{cor}

\begin{prf*}
  The inequality $\sup\{\Gfcd M \mid M\text{ is an $A$-module}\} \le
  \Ggldim$ holds by \thmcite[3.3]{CET-21a}, so the assertion is a
  special case of \prpref{supGfcd}.
\end{prf*}

This gives new examples of rings, including non-coherent rings, with
the property that strongly cotorsion modules and Gorenstein injective
modules coincide.

\begin{exa}
  Over any ring of finite global dimension the strongly cotorsion
  modules, injective modules, and Gorenstein injective modules
  coincide. A result of Enochs, Estrada, and Iacob \cite[Theorem
  3.2]{EEI-08} yields less trivial examples: the ring of dual numbers
  over any ring of finite Gorenstein global dimension. A concrete
  example of a non-coherent ring of finite global dimension is
  provided by Estrada, Iacob, and Yeomans \cite[Section~4(1)]{EIY-17}.
\end{exa}

Yoshizawa \cite{TYs12} in fact shows that among (commutative)
Cohen-Macaulay complete local rings, those that are Gorenstein are
characterized by the property that strongly cotorsion modules and
Gorenstein injective modules are the same. More generally, Iacob
\cite{AIc17} shows that this property characterizes Iwanaga-Gorenstein
rings among noetherian rings $A$ with $\FPD<\infty$ and
$\textrm{FPD}(\Aop)<\infty$. In view of \prpref{GI_is_SC}, the
characterizing property is really that strongly cotorsion modules are
Gorenstein injective. Further, a noetherian ring $A$ is
Iwanaga-Gorenstein if and only if $\Ggldim$ is finite, so the next
result can be interpreted as removing the noetherian assumption from
\thmcite[10]{AIc17}:

\begin{thm}\label{thm:Ggl_SC_GI}
  The following conditions are equivalent.
  \begin{eqc}
  \item $\Ggldim<\infty$.
  \item $\FPD<\infty$ and every strongly cotorsion $A$-module is
    Gorenstein injective.
  \end{eqc}
\end{thm}

\begin{prf*}
  The implication \proofofimp[]{i}{ii} holds by
  \corref{Ggldim_fin_SCGI} and \thmcite[2.28]{HHl04a} which shows that
  the assumption on $A$ implies that $\FPD$ is finite. For
  \proofofimp[]{ii}{i}, set $n=\FPD$. By a result of Jensen
  \prpcite[6]{CUJ70}, every flat $A$-module, and hence every module of
  finite flat dimension, has projective dimension at most $n$.  It
  suffices to show that $\Gid{N}\le n$ holds for every $A$-module
  $N$. Let $N$ be an $A$-module and consider an exact sequence,
  \begin{equation*} 0 \lra N \lra E_0\lra \cdots \lra E_{-n+1} \lra D
    \lra 0 \:,
  \end{equation*}
  with each $E_i$ injective.  For every $A$-module $M$ with
  $\fd{M}\le n$, dimension shifting along this sequence yields
  $\Ext{1}{M}{D}\cong \Ext{n+1}{M}{N}=0$. Therefore, $D$ is a strongly
  cotorsion $A$-module. By assumption $D$ is thus Gorenstein
  injective, whence $\Ggldim$ is finite.
\end{prf*}

For ease of comparison to the literature, we include:
\begin{cor}
  If $\FPD$ is finite, then the following conditions are equivalent.
  \begin{eqc}
  \item $\Ggldim<\infty$.
  \item An $A$-module is strongly cotorsion if and only if it is
    Gorenstein injective.
  \end{eqc}
\end{cor}

\begin{rmk}
  In the literature it is often remarked that Ding-Chen rings have
  properties that generalize those of Iwanaga-Gorenstein rings to a
  non-noetherian setting; see for example the abstract of
  \cite{JGl10}. Ding and Chen \thmcite[7]{NDnJCh96} show that a
  Ding-Chen ring has finite Gorenstein weak global dimension. In light
  of \rmkref{Gil}, however, \thmref{Ggl_SC_GI} says that a Ding-Chen
  ring $A$ has finite Gorenstein global dimension if and only if
  $\FPD<\infty$ holds. It follows from \prpcite[6]{CUJ70} that a von
  Neumann regular ring $A$ of infinite global dimension is a a
  Ding-Chen ring with $\FPD=\infty$. By a result of Pierce
  \cite[Corollary 5.2]{RSP67}, a free boolean ring with
  $\aleph_{\omega}$ generators, for an infinite cardinal $\omega$, is
  a von Neumann regular ring of infinite global dimension. We remark
  that this observation was essentially already made by Wang
  \cite[Example 3.3]{JWn17}.
\end{rmk}

Here is another consequence of \prpref{supGfcd}; without assumptions
on the ring it does not readily compare to \corref{Ggldim_fin_SCGI}.

\begin{cor}
  \label{cor:Gwgldim_fin_SCGI}
  If $\Gwgldim$ is finite, then an $A$-module is strongly cotorsion if
  and only if it is Gorenstein injective.
\end{cor}

\begin{prf*}
  The inequality $\sup\{\Gfcd M \mid M\text{ is an $A$-module}\} \le
  \Gwgldim$ holds by \thmcite[5.7]{CELTWY-21}, so the assertion is a
  special case of \prpref{supGfcd}.
\end{prf*}

\begin{thm}
  \label{thm:Gfcd_SC_GI}
  If $A$ is left or right coherent, then the next conditions are
  equivalent.
  \begin{eqc}
  \item $\Gwgldim<\infty$.
  \item $\FFD<\infty$ and every strongly cotorsion $A$-module is
    Gorenstein injective.
  \end{eqc}
\end{thm}

\begin{prf*}
  Condition $(i)$ implies $(ii)$ by \corref{Gwgldim_fin_SCGI} and
  \thmcite[3.24]{HHl04a} as the Gorenstein weak global dimension is
  symmetric by \corcite[2.5]{CET-21a}. For the converse, assume that
  $A$ is left coherent and set $n = \FFD$. Let $M$ be an $\Aop$-module
  with a flat resolution $F$, and let $N=\Co[n]{F}$ be the
  $n^\mathrm{th}$ syzygy in this resolution. By \thmcite[3.6]{HHl04a}
  the module $N$ is Gorenstein flat if and only if the dual
  $N^+=\Hom[\mathbb{Z}]{M}{\mathbb{Q}/\mathbb{Z}}$ is a Gorenstein
  injective $A$-module. Let $L$ be an $A$-module of finite flat
  dimension. The modules $F_i^+$ are injective, so dimension shifting
  along
  \begin{equation*}
    0 \lra M^+ \lra (F_0)^+ \lra \cdots \lra (F_{n-1})^+ \lra N^+\lra 0
  \end{equation*}
  yields $\Ext{1}{L}{N^+}\cong \Ext{n+1}{L}{M^+}$. As $\fd L\le n$ and
  $M^+$ is cotorsion, dimension shifting along a flat resolution of
  $L$ yields $\Ext{n+1}{L}{M^+}=0$. Thus $N^+$ is strongly cotorsion
  and hence Gorenstein injective. This shows that every $\Aop$-module
  has Gorenstein flat dimension at most $n$, and
  $\Gwgldim = \Gwgldim[\Aop]$ holds by \corcite[2.5]{CET-21a}. If $A$
  is instead right coherent, then the same argument applies with $A$ and
  $\Aop$ interchanged.
\end{prf*}

Again for ease of comparison to the literature, we include:

\begin{cor}
  If $A$ is left or right coherent and $\FFD$ is finite, then the
  following conditions are equivalent.
  \begin{eqc}
  \item $\Gwgldim<\infty$.
  \item An $A$-module is strongly cotorsion if and only if it is
    Gorenstein injective.
  \end{eqc}
\end{cor}

For modules over a right coherent ring, the Gorenstein flat dimension
agrees with the Gorenstein flat-cotorsion dimension, see
\corcite[5.8]{CELTWY-21}. In view of \thmref{Gfcd_SC_GI} it is thus
natural to ask if \prpref{supGfcd} has a converse:

\begin{qst}
  If $\FFD$ is finite and an $A$-module is strongly
  cotorsion if and only if it is Gorenstein injective, is
  $\sup\{\Gfcd M \mid M\text{ is an $A$-module}\}$ finite?
\end{qst}

\begin{rmk}
  We remark that finiteness of
  $\sup\{\Gfcd M \mid M\text{ is an $A$-module}\}$ implies that
  another two quantities are finite, namely $\sfli$, as already
  noticed in the proof of \prpref{supGfcd}, and
  $\sup\{\id{M} \mid M\text{ is a flat-cotorsion $A$-module}\}$. If
  one adds to the assumption $\FFD < \infty$ the assumption that these
  two quantities are finite, then it seems feasible to prove that
  $\sup\{\Gfcd M \mid M\text{ is an $A$-module}\}$ is finite, but at
   this time we are not convinced that this is the best one can do.
\end{rmk}

%%%%%%%%%%%%%%%%%%%%%%%%%%%%%%%%%%%%%%%%%%%%%%%%%%%%%%%%%%%%%%%%%%%%%%
\section{Group rings}
\label{sec_app}

\noindent
We take the opportunity to include a note on how the Gorenstein
flat-cotorsion dimension behaves along certain ring
homomorphisms. This leads to an answer to a question of Y.~Xiang
\cite{YXn21} on the Gorenstein weak global dimension of group rings.

First make the following observation: Let $A\to B$ be a flat ring
homomorphism. If $C$ is a cotorsion $B$-module, then it is also
cotorsion as an $A$-module. To see this, let $F$ be a flat $A$-module
and note that standard Hom-tensor adjunction yields
$\Ext[A]{1}{F}{C}\cong \Ext[A]{1}{F}{\Hom[B]{B}{C}}\cong
\Ext[B]{1}{B\otimes_A F}{C}=0$.

\begin{prp}
  \label{prp:freeflat}
  Let $A\to B$ be a homomorphism of rings and $M$ a $B$-module. If $B$
  is free as an $A$-module and $\Gfcd[B]{M}<\infty$, then
  $\Gfcd[A]{M}\ge \Gfcd[B]{M}$.
\end{prp}

\begin{prf*}
  Set $n=\Gfcd[B]{M}$. By \thmcite[4.5]{CELTWY-21}, one thus has
  $\Ext[B]{n}{M}{C}\not=0$ for some flat-cotorsion $B$-module $C$. By
  the earlier remark, $C$ is also flat-cotorsion viewed as an
  $A$-module.  It now follows that
$$\Ext[A]{n}{M}{C}\cong \Ext[A]{n}{B\otimes_B M}{C}\cong \Ext[B]{n}{M}{\Hom[A]{B}{C}},$$
which contains $\Ext[B]{n}{M}{C}$ as a direct summand and hence is
non-zero. Thus $n\le \Gfcd[A]{M}$.
\end{prf*}

As an application, the following result removes the right coherence
assumption in \prpcite[4.9]{YXn21}, thus affirmatively answering
\cite[Question 4.11]{YXn21}.

\begin{thm}
  \label{thm:grouprings}
  Let $k$ be a field, $G$ a group, and $H$ a subgroup of $G$ of finite
  index. If $\Gwgldim[{k[G]}]$ is finite, then it is equal to
  $\Gwgldim[{k[H]}]$.
\end{thm}

\begin{prf*}
  The inequality $\Gwgldim[{k[H]}] \le \Gwgldim[{k[G]}]$ holds by
  \prpcite[4.3(2)]{YXn21}. Assume that $\Gwgldim[{k[G]}]=n<\infty$.
  By \prpcite[4.2]{YXn21}, one has $\Gwgldim[{k[G]}]=\Gfd[{k[G]}]{k}$,
  which in turn equals $\Gfcd[{k[G]}] k$ since Gorenstein
  flat-cotorsion and Gorenstein flat dimensions agree when the latter
  is finite by \thmcite[5.7]{CELTWY-21}. Since $k[G]$ is free as a
  $k[H]$-module, \prpref{freeflat} yields
  $\Gfcd[{k[G]}] k \le \Gfcd[{k[H]}] k$, and so we obtain the other
  inequality: $\Gwgldim[{k[H]}] \ge \Gwgldim[{k[G]}]$.
\end{prf*}

% \section*{Acknowledgment}

% \noindent
% We thank the referee for pertinent suggestions that improved the
% exposition.

% \bibliographystyle{amsplain} \bibliography{../+references}
%\bibliographystyle{../amsplain-nodash} \bibliography{../+references}

\def\soft#1{\leavevmode\setbox0=\hbox{h}\dimen7=\ht0\advance \dimen7
  by-1ex\relax\if t#1\relax\rlap{\raise.6\dimen7
  \hbox{\kern.3ex\char'47}}#1\relax\else\if T#1\relax
  \rlap{\raise.5\dimen7\hbox{\kern1.3ex\char'47}}#1\relax \else\if
  d#1\relax\rlap{\raise.5\dimen7\hbox{\kern.9ex \char'47}}#1\relax\else\if
  D#1\relax\rlap{\raise.5\dimen7 \hbox{\kern1.4ex\char'47}}#1\relax\else\if
  l#1\relax \rlap{\raise.5\dimen7\hbox{\kern.4ex\char'47}}#1\relax \else\if
  L#1\relax\rlap{\raise.5\dimen7\hbox{\kern.7ex
  \char'47}}#1\relax\else\message{accent \string\soft \space #1 not
  defined!}#1\relax\fi\fi\fi\fi\fi\fi}
  \providecommand{\MR}[1]{\mbox{\href{http://www.ams.org/mathscinet-getitem?mr=#1}{#1}}}
  \renewcommand{\MR}[1]{\mbox{\href{http://www.ams.org/mathscinet-getitem?mr=#1}{#1}}}
  \providecommand{\arxiv}[2][AC]{\mbox{\href{http://arxiv.org/abs/#2}{\sf
  arXiv:#2 [math.#1]}}} \def\cprime{$'$}
\providecommand{\bysame}{\leavevmode\hbox to3em{\hrulefill}\thinspace}
\providecommand{\MR}{\relax\ifhmode\unskip\space\fi MR }
\providecommand{\MRhref}[2]{%
  \href{http://www.ams.org/mathscinet-getitem?mr=#1}{#2}
}
\providecommand{\href}[2]{#2}

\end{document}